\documentclass[a4paper, 11pt]{article}

\usepackage{amsmath,amsfonts}
\usepackage{amsthm}
\usepackage{amssymb}
\usepackage{bm}

\usepackage{mathtools}
\mathtoolsset{showonlyrefs=true}

\usepackage{tikz-cd}
\usepackage{graphicx}
\usepackage{color}
\usepackage{xcolor}
\usepackage{listings}
\usepackage{latexsym}

\usepackage{hyperref}

\hypersetup{
setpagesize=false,
bookmarksnumbered=true,%
bookmarksopen=true,%
colorlinks=false,
linkcolor=red,
citecolor=blue,
urlcolor=orange,
}

\usepackage{titling}

\usepackage{tcolorbox}
\tcbuselibrary{breakable, skins, theorems}

\newif\ifdraft
\draftfalse 

\ifdraft
  \definecolor{burgundy}{rgb}{0.5, 0.0, 0.13}

  \newtcolorbox{mydfnbox}{empty, breakable,
  underlay={
  \fill[blue!10] (frame.north west) rectangle (frame.south east);
  \draw[line width=2pt] ([xshift=1pt]frame.north west)--([xshift=1pt]frame.south west);
  }}

  \newtcolorbox{mythmbox}{empty, breakable,
  underlay={
  \fill[red!10!white] (frame.north west) rectangle (frame.south east);
  \draw[line width=2pt] ([xshift=1pt]frame.north west)--([xshift=1pt]frame.south west);
  }}

  \newtcolorbox{mylemmabox}{empty, breakable,
  underlay={
  \fill[green!10] (frame.north west) rectangle (frame.south east);
  \draw[line width=2pt] ([xshift=1pt]frame.north west)--([xshift=1pt]frame.south west);
  }}

  \definecolor{charcoal}{rgb}{0.21, 0.27, 0.31}
  \newtcolorbox{myexamplebox}{empty, breakable,
  underlay={
  \fill[charcoal!20] (frame.north west) rectangle (frame.south east);
  \draw[line width=2pt] ([xshift=1pt]frame.north west)--([xshift=1pt]frame.south west);
  }}

  \newtcolorbox{myquebox}{empty, breakable,
  underlay={
  \fill[yellow!10] (frame.north west) rectangle (frame.south east);
  \draw[line width=2pt] ([xshift=1pt]frame.north west)--([xshift=1pt]frame.south west);
  }}
 
  \newtheorem{mydfn}{Definition}[section]
  \newtheorem{myprop}{Proposition}[section]
  \newtheorem{mythm}{Theorem}[section]
  \newtheorem{mylemma}{Lemma}[section]
  \newtheorem{myconj}{Conjecture}[section]
  \newtheorem{myrem}{Remark}[section]
  \newtheorem{mycor}{Corollary}[section]
  \newtheorem{myque}{Question}[section]
  \newtheorem{myexample}{Example}[section]

  \newenvironment{dfn}{\begin{mydfnbox}\begin{mydfn}}{\end{mydfn}\end{mydfnbox}}
  \newenvironment{prop}{\begin{mythmbox}\begin{myprop}}{\end{myprop}\end{mythmbox}}
  \newenvironment{thm}{\begin{mythmbox}\begin{mythm}}{\end{mythm}\end{mythmbox}}
  \newenvironment{lemma}{\begin{mylemmabox}\begin{mylemma}}{\end{mylemma}\end{mylemmabox}}
  \newenvironment{conj}{\begin{myquebox}\begin{myconj}}{\end{myconj}\end{myquebox}}
  \newenvironment{rem}{\begin{myexamplebox}\begin{myrem}}{\end{myrem}\end{myexamplebox}}
  \newenvironment{cor}{\begin{mythmbox}\begin{mycor}}{\end{mycor}\end{mythmbox}}
  \newenvironment{que}{\begin{myquebox}\begin{myque}}{\end{myque}\end{myquebox}}
  \newenvironment{example}{\begin{myexamplebox}\begin{myexample}}{\end{myexample}\end{myexamplebox}}
\else
  \theoremstyle{definition}
  \newtheorem{dfn}{Definition}[section]
  \newtheorem{prop}{Proposition}[section]
  \newtheorem{thm}{Theorem}[section]
  \newtheorem{lemma}{Lemma}[section]
  
  \newtheorem{rem}{Remark}[section]
  \newtheorem{cor}{Corollary}[section]
  
  \newtheorem{example}{Example}[section]
\fi

\begin{document}

\title{Takai duality and crossed product hosts\\ in C*-actions}
\author{
  \Large Yusuke Nakae\thanks{
    Graduate School of Mathematical Sciences, The University of Tokyo, Tokyo, Japan.\\
    Email: nakae@ms.u-tokyo.ac.jp
    }
}
\date{}

\maketitle

\vspace{-3em}

\begin{abstract}
Crossed product algebras are fundamental in the study of C*-algebras, traditionally under the assumption of continuity of group actions. Recent work by Grundling and Neeb introduced the crossed product host, an analog of the crossed product for a singular action. 

In this paper, we investigate the structure of the crossed product host and its relation to the conventional crossed product. We examine the validity of Takai-type duality in this setting and establish connections via the Landstad algebra. Additionally, we provide a necessary and sufficient condition for the existence of ground states, and show that for amenable groups, if the full and reduced crossed product hosts exist, then they coincide.

\noindent \textbf{Mathematics Subject Classification:} 22D25, 46L55, 47L65, 81R15.
\end{abstract}

\tableofcontents

\vspace{2em}
\renewcommand{\include}[1]{}
\renewcommand\documentclass[2][]{}

\section{Introduction}

Crossed product algebras have played a crucial role in the study of C*-algebras \cite{MR2288954} and von Neumann algebras. When we consider crossed products in the context of C*-algebras, they have traditionally been assumed that the group actions are continuous. However, in applications, it would be more useful to handle cases where the group actions are not necessarily continuous.

Recently, Grundling and Neeb, among others, have developed an extended approach to studying C*-dynamical systems with such singular actions from a representation-theoretic perspective \cite{MR3177335, MR4157367}. They refer to such systems as C*-actions. The corresponding analog of the crossed product algebra in their framework is known as the crossed product host. Unlike in the regular case, some desirable properties may no longer hold, yet many analogous structures and results still persist.

One of the most fundamental results in the crossed product theory is the Takai duality for abelian groups \cite{MR365160, MR989764} and its extension to non-abelian groups, known as the Imai-Takai duality \cite{MR500719}. A natural question arises: Does a similar duality hold in the setting of crossed product hosts? Additionally, we aim to explore methods for interpreting the somewhat elusive crossed product host in terms of more regular structures.

In Section 2, we introduce the definition of the crossed product host along with several related notions. Section 3 is devoted to our main results concerning duality; we investigate to what extent a Takai-type duality holds in this framework. Finally, we utilize the Landstad conditions to express the crossed product host in terms of a conventional crossed product algebra and examine comparisons between the original C*-action and the derived C*-dynamical system. Additionally, using these results, we establish the existence of ground states within the crossed product host framework. Furthermore, we consider the reduced version of the crossed product host, analogous to the standard crossed product algebra, and clarify its relationship with the full version. Section 4 concludes with a summary of our findings and directions for future research.
\section{Preliminaries}

We use a non-degenerate homomorphism in the following sense, from \cite{MR989764} on page 2. 
\begin{dfn}[A non-degenerate homomorphism]
  Let $M(\mathcal{B})$ denote the multiplier algebra of a C*-algebra $\mathcal{B}$. We say a homomorphism $\phi:\mathcal{A} \to M(\mathcal{B})$ is \textbf{non-degenerate} if there is an approximate identity $(u_{\lambda})_{\lambda}$ for $\mathcal{A}$ such that $\phi(u_{\lambda}) \to 1$ strictly in $M(\mathcal{B})$.
\end{dfn}

\begin{rem}
  In the above situation, there is the unique non-degenerate strictly continuous homomorphism from $M(\mathcal{A})$ to $M(\mathcal{B})$, denoted as $\widetilde{\phi}$.
\end{rem}

For a complex Hilbert space $\mathcal{H}$ and a C*-algebra $\mathcal{A}$, we write $\mathrm{Rep}(\mathcal{A},\mathcal{H})$ for the set of non-degenerate representations of $\mathcal{A}$ on $\mathcal{H}$. For a topological group $G$, we write $\mathrm{Rep}(G,\mathcal{H})$ for the set of strongly continuous unitary representations of $G$ on $\mathcal{H}$.

Next, we define a C*-action and a host algebra. They are introduced in \cite{MR3177335}. 

\begin{dfn}[A C*-action]
  Let $\mathcal{A}$ be a C*-algebra, $G$ be a locally compact group, $\alpha$ be an action of $G$ on $\mathcal{A}$.
  We call a triple $(\mathcal{A}, G, \alpha)$ is \textbf{C*-action}.
  (There is no continuity condition for $\alpha$.) 

  When the action $\alpha$ is continuous, the triple $(\mathcal{A}, G, \alpha)$ is called a \textbf{C*-dynamical system}.
\end{dfn}

\begin{dfn}[A host algebra]
  Let $G$ be a topological group. A pair $(\mathcal{L}, \eta)$ is called a \textbf{host algebra} where $\mathcal{L}$ is a C*-algebra and $\eta:G\to U(M(\mathcal{L}))$ is a group homomorphism satisfying the following two conditions.

  (H1) For a non-degenerate representation $(\pi, \mathcal{H})$ of $\mathcal{L}$, the representation $(\widetilde{\pi}\circ \eta, \mathcal{H})$ of $G$ is (strongly) continuous.

  (H2) For each complex Hilbert space $\mathcal{H}$, the map
  \begin{align*}
    \eta_{*}: \mathrm{Rep}(\mathcal{L},\mathcal{H}) \to \mathrm{Rep}(G, \mathcal{H}), \pi\mapsto \widetilde{\pi}\circ \eta
  \end{align*}
  is injective.

  Moreover, if it satisfies the condition below, it is called a \textbf{full} host algebra.

  (H3) The map $\eta_{*}$ is surjective.
\end{dfn}

For a locally compact group $G$, we only consider the full host algebra $\mathcal{L}=C^*(G)$. In this case, we take the canonical map $\eta: G \to U(M(C^*(G)))$. 

Next we introduce the notion of a crossed product host, from \cite[Definition 3.2]{MR4157367} on page 8. This is also introduced in \cite{MR3177335}.
\begin{dfn}[A crossed product host]
  Let $G$ be a topological group, and let $(\mathcal{L}, \eta)$ be a host algebra for $G$ and $(\mathcal{A},G,\alpha)$ be a C*-action. We call a triple $(\mathcal{C}, \eta_{\mathcal{A}},\eta_{\mathcal{L}})$ a \textbf{crossed product host} for $(\alpha, \mathcal{L})$ if it satisfies the following four conditions.
  
  (CP1) The maps $\eta_{\mathcal{A}}:\mathcal{A} \to M(\mathcal{C})$ and $\eta_{\mathcal{L}}:\mathcal{L} \to M(\mathcal{C})$ are morphisms of C*-algebras.

  (CP2) The map $\eta_{\mathcal{L}}$ is non-degenerate, i.e., $\eta_{\mathcal{L}}(\mathcal{L})\mathcal{C}$ is dense in $\mathcal{C}$.

  (CP3) The multiplier extension $\widetilde{\eta}_{\mathcal{L}}:M(\mathcal{L}) \to M(\mathcal{C})$ satisfies in $M(\mathcal{C})$ the relations
  \begin{align*}
    \widetilde{\eta}_{\mathcal{L}}(\eta(g))\eta_{\mathcal{A}}(A)\widetilde{\eta}_{\mathcal{L}}(\eta(g))^* = \eta_{\mathcal{A}}(\alpha_{g}(A)) \quad \text{for all} \quad A\in \mathcal{A}, g\in G.
  \end{align*}

  (CP4) It holds $\eta_{\mathcal{A}}(\mathcal{A})\eta_{\mathcal{L}}(\mathcal{L}) \subseteq \mathcal{C}$ and $\mathcal{C}$ is generated by this set as a C*-algebra.

  We say the crossed product host for $(\alpha,\mathcal{L})$ is \textbf{full} if it also satisfies

  (CP5) For every covariant representation $(\pi, U)$ of $\mathcal{A}$ on $\mathcal{H}$ for which $U$ is an $\mathcal{L}$-representation of $G$, there exists a unique representation $\rho:\mathcal{C} \to B(\mathcal{H})$ with
  \begin{align*}
    \rho(\eta_{\mathcal{A}}(A)\eta_{\mathcal{L}}(L))=\pi(A)U_{\mathcal{L}}(L), \text{ for } A\in \mathcal{A}, L\in \mathcal{L}.
  \end{align*}
\end{dfn}

If a C*-action $(\mathcal{A}, G, \alpha)$ has a full crossed product host for $(\alpha, C^*(G))$, we write $\mathrm{fcph}(\mathcal{A}, G, \alpha)$ for it.

\begin{dfn}[A cross representation]\label{def: cross rep}
    Let $(\mathcal{L},\eta)$ be a host algebra for a locally compact group $G$. We say a covariant representation $(\pi, U)$ of a C*-action $(\mathcal{A}, G, \alpha)$ on $\mathcal{H}$ is \textbf{cross} when it satisfies the following equivalent conditions.

    (i) The triple $(\mathcal{C}:=C^*(\pi(\mathcal{A})U_{\mathcal{L}}(\mathcal{L})), \eta_{\mathcal{A}}, \eta_{\mathcal{L}})$ is a crossed product host.
    
    (ii) It holds $\pi(\mathcal{A})U_\mathcal{L}(\mathcal{L}) \subseteq U_\mathcal{L}(\mathcal{L})B(\mathcal{H})$.

    (iii) For every approximate identity $(E_{j})_{j\in J}$ of $\mathcal{L}$ we have
    \begin{align*}
      \|U_{\mathcal{L}}(E_{j})\pi(A)U_{\mathcal{L}}(L) - \pi(A)U_{\mathcal{L}}(L)\| \to 0 \quad \text{for} \quad A\in \mathcal{A}, L\in \mathcal{L}. 
    \end{align*}

    (iv) There exists an approximate identity $(E_{j})_{j\in J}$ of $\mathcal{L}$ such that
    \begin{align*}
      \|U_{\mathcal{L}}(E_{j})\pi(A)U_{\mathcal{L}}(L) - \pi(A)U_{\mathcal{L}}(L)\| \to 0 \quad \text{for} \quad A\in \mathcal{A}, L\in \mathcal{L}. 
    \end{align*}
\end{dfn}

We now define a coaction, adopted from the definition of \cite[Appendix A.3, page 127]{MR2203930}. The notation $\mathcal{A}\otimes C^*(G)$ means the minimal tensor product of $\mathcal{A}$ and $C^*(G)$.

\begin{dfn}[A coaction]
Let $\mathcal{A}$ be a C*-algebra and $G$ a locally compact group. An injective, non-degenerate $*$-homomorphism 
$$
\delta: \mathcal{A} \to M(\mathcal{A} \otimes C^*(G))
$$ 
is called a \textbf{coaction} if it satisfies:
\begin{enumerate}
  \item[(i)] $\delta(\mathcal{A})(1 \otimes C^*(G)) \subseteq \mathcal{A} \otimes C^*(G)$,
  \item[(ii)] $(\delta \otimes \mathrm{id}_G)\circ \delta = (\mathrm{id}_{\mathcal{A}} \otimes \delta_G)\circ \delta$,
\end{enumerate}
where $\delta_G: C^*(G) \to M(C^*(G) \otimes C^*(G))$ denotes the canonical comultiplication on $C^*(G)$.

Furthermore, the coaction $\delta$ is said to be \textbf{non-degenerate} if, in addition, it satisfies:
\begin{enumerate}
  \item[(iii)] $\overline{\delta(\mathcal{A})(1\otimes C^*(G))} = \mathcal{A} \otimes C^*(G)$.
\end{enumerate}
\end{dfn}

Next, in order to describe duality, we define the notion of a maximal coaction.

\begin{dfn}[A maximal coaction]
A coaction $\delta$ of $G$ on $\mathcal{A}$ is said to be \textbf{maximal} if
\begin{align*}
  \mathcal{A} \rtimes_{\delta} G \rtimes_{\hat{\delta}} G \simeq \mathcal{A} \otimes K(L^2(G)),
\end{align*}
where, $K(L^2(G))$ is the set of all compact operators on the Hilbert space $L^2(G)$.
\end{dfn}
\section{Main Result}

\subsection{An iterated crossed product host}

In this subsection, we assume that $G$ is a locally compact abelian group. We consider a C*-action $(\mathcal{A}, G, \alpha)$. and we assume that it has the full crossed product host, denoted as $\mathrm{fcph}(\mathcal{A}, G, \alpha)$, as remarked in the previous section. Firstly, we show the universal property of the full crossed product host.

\begin{prop}\label{prop: universal property of fcph}
  Let $\phi$ be a non-degenerate homomorphism of $\mathcal{A}$ into the multiplier algebra $M(\mathcal{C})$ of a C*-algebra $\mathcal{C}$, and $U$ be a strictly continuous homomorphism of $G$ into $UM(\mathcal{C})$, the unitaries of $M(\mathcal{C})$, such that $\phi(\alpha_{s}(A))=U(s)\phi(A)U(s)^*$ for all $A\in \mathcal{A}, s\in G$. Then there is a non-degenerate homomorphism $\rho$ of $\mathrm{fcph}(\mathcal{A},G,\alpha)$ into $M(\mathcal{C})$ such that $\widetilde{\rho} \circ \eta_{\mathcal{A}} = \phi$ and $\widetilde{\rho}\circ \eta_{\mathcal{L}}\circ \eta = U$.
\[\begin{tikzcd}
	{\mathcal{A}} & {M(\mathrm{fcph}(\mathcal{A},G,\alpha))} & {M(\mathcal{L})} & G \\
	& {M(\mathcal{C})}
	\arrow["{\eta_\mathcal{A}}", from=1-1, to=1-2]
	\arrow["\phi"', from=1-1, to=2-2]
	\arrow["{\exists\widetilde{\rho}}", dashed, from=1-2, to=2-2]
	\arrow["{\widetilde{\eta}_\mathcal{L}}"', from=1-3, to=1-2]
	\arrow["\eta"', from=1-4, to=1-3]
	\arrow["U", from=1-4, to=2-2]
\end{tikzcd}\]
\end{prop}
\begin{proof}
  Let $(\pi, \mathcal{H})$ be a faithful representation of $\mathcal{C}$. Then $(\pi\circ \phi, \pi \circ U)$ is a covariant representation of $(\mathcal{A},G,\alpha)$. Because of the full crossed product property, there exists a representation $(\rho_{\mathcal{H}},\mathcal{H})$ of $\mathrm{fcph}(\mathcal{A},G,\alpha)$ for which the following holds.
  \begin{align*}
    \rho_{\mathcal{H}}\circ \eta_\mathcal{A} = \pi\circ \phi, \\
    \rho_\mathcal{H}\circ \eta_G = \pi\circ U.
  \end{align*}
  Therefore, $\rho:=\pi^{-1}|_{\pi(\mathcal{C})}\circ \rho_{\mathcal{H}}$ satisfies the conditions. Non-degeneracy is easily verified.
\end{proof}

The dual group $\hat{G}$ acts trivially on $\mathcal{A}$. We write this trivial action for $\iota$. Then we can construct the crossed product $\mathcal{A}\rtimes_{\iota} \hat{G}$ (this is isomorphic to $C_{0}(G,\mathcal{A})$). An action $\beta$ is induced by $\alpha$ on $\mathcal{A}\rtimes_{\iota} \hat{G}$ with
\begin{align*}
  (\beta_{g}(f))(\gamma):=\langle g, \gamma\rangle\alpha_{g}(f(\gamma)), \text{ for } f\in C_{c}(\hat{G}, \mathcal{A}) \subset \mathcal{A}\rtimes_{\iota} \hat{G}, g\in G, \gamma\in \hat{G}. 
\end{align*}

In this way, $(\mathcal{A} \rtimes_{\iota} \hat{G}, G, \beta)$ forms a C*-action.

Next, we show that a cross representation in $(\mathcal{A},G,\alpha)$ is inherited as a cross representation in $(\mathcal{A} \rtimes_{\iota} \hat{G}, G, \beta)$.

\begin{prop}
 Let $(\mathcal{A}, G, \alpha)$ be a C*-action. If a covariant representation $(\pi, U)$ of $(\mathcal{A}, G, \alpha)$ is cross, then the covariant representation $(\pi \rtimes V, U)$ of C*-action $(\mathcal{A} \rtimes_{\iota} \hat{G}, G, \beta)$ on $\mathcal{H}$ is also cross, where $V$ is a unitary representation of $\hat{G}$. 
\end{prop}
\begin{proof}
 By the definition of a cross representation (Definition \ref{def: cross rep} (ii)), it is sufficient to show that
 \begin{align*}
     (\pi\times V)(\mathcal{A}\rtimes \hat{G})U_{\mathcal{L}'}(\mathcal{L}') \subseteq U_{\mathcal{L}'}(\mathcal{L}')B(\mathcal{H}).
 \end{align*}
 Moreover, it is sufficient to show that
 \begin{align*}
     (\pi\times V)(C_c(\hat{G},\mathcal{A}))U_{\mathcal{L}'}(C_c(G)) \subseteq U_{\mathcal{L}'}(\mathcal{L}')B(\mathcal{H}).
 \end{align*}
 Moreover, it is sufficient to show that for all $A\in \mathcal{A}, \varphi\in C_c(\hat{G}), f\in C_c(G)$
 \begin{align*}
     (\pi\times V)(A\otimes \varphi)U_{\mathcal{L}'}(f) \in U_{\mathcal{L}'}(\mathcal{L}')B(\mathcal{H}).
 \end{align*}
 Then,
 \begin{align*}
     (\pi\times V)(A\otimes \varphi)U_{\mathcal{L}'}(f)&=\pi(A)V(\varphi)U(f)\\
     &=\pi(A)\int \varphi(\gamma)V_\gamma d\gamma \int f(g) U(g)dg\\
     &=\pi(A)\int\int \langle g, \gamma\rangle^{-1}f(g) U_{g} \varphi(\gamma)V_\gamma dg d\gamma\\
     &=\int \pi(A)U(f_{-\gamma}) \varphi(\gamma)V_{\gamma} d \gamma,
 \end{align*}
 where, we have set $f_{-\gamma}(g):=\langle g, \gamma\rangle^{-1}f(g)$. Then, since $\pi(A)U(f_{-\gamma}) \varphi(\gamma)V_{\gamma} \in U(\mathcal{L})B(\mathcal{H})$, so the total is also in $U(\mathcal{L})B(\mathcal{H})$.
\end{proof}

\begin{cor}
    If a C*-action $(\mathcal{A}, G, \alpha)$ has the full crossed product host, then the C*-action $(\mathcal{A} \rtimes_{\iota} \hat{G}, G, \beta)$ has also the full crossed product host.
\end{cor}
\begin{proof}
    Since any non-degenerate representation of $\mathcal{A} \rtimes_{\iota} \hat{G}$ arises from a covariant representation of $(\mathcal{A}, G, \alpha)$, their covariant representations are all cross. Hence, it has the full crossed product.
\end{proof}

It is well known that in the regular case, i.e., $G$ is an abelian and locally compact group and we consider a C*-dynamical system $(\mathcal{A}, G, \alpha)$, then the following isomorphism holds:
\begin{align*}
  \mathcal{A} \rtimes_{\alpha} G \rtimes_{\hat{\alpha}} \hat{G}\simeq \mathcal{A} \rtimes_{\iota} \hat{G} \rtimes_{\beta} G.
\end{align*}
A natural question arises: what happens in the singular case?

In fact, we have established the following proposition.

\begin{thm}(Duality of an iterated crossed product host)
  Let $G$ be a locally compact abelian group, $(\mathcal{A}, G, \alpha)$ be a C*-action. We assume that this has the full crossed product host $\mathrm{fcph}(\mathcal{A}, G, \alpha)$. Then, we have
  \begin{align*}
    \mathrm{fcph}(\mathcal{A},G,\alpha) \rtimes_{\hat{\alpha}} \hat{G} \simeq \mathrm{fcph}(\mathcal{A}\rtimes_{\iota} \hat{G},G,\beta).
  \end{align*}
\end{thm}

The proof below is inspired by the arguments in \cite{MR2288954}.

\begin{proof}
We construct a mapping and show that it is a homomorphism and bijective. For simplicity, we write \(\mathcal{C}\) and \(\mathcal{D}\) for 
\begin{align*}
\mathrm{fcph}(A,G,\alpha) \rtimes_{\hat{\alpha}} \hat{G} \quad \text{and} \quad \mathrm{fcph}(\mathcal{A}\rtimes_{\iota} \hat{G},G,\beta)
\end{align*}
respectively.

\medskip

\textbf{(Construction of the map $\Phi_{\mathcal{D}}$)}

We construct the mapping by repeatedly applying the universal property. First, a straightforward calculation shows that
\begin{align*}
C_{b}(\hat{G}, M(\mathrm{fcph}(\mathcal{A}, G, \alpha))) \subset M(\mathcal{C}),
\end{align*}
so we regard the left-hand side as a subalgebra and construct a mapping whose codomain consists of functions on \(\hat{G}\) rather than \(M(\mathcal{C})\).

First, define
\begin{align*}
\Phi_{\mathcal{A}}: \mathcal{A} \to C_{b}(\hat{G}, M(\mathrm{fcph}(\mathcal{A}, G, \alpha)))
\end{align*}
by
\begin{align*}
\Phi_{\mathcal{A}}(A) := \Bigl(\tau \mapsto \eta_{\mathcal{A}}(A)\widetilde{\eta}_{\mathcal{L}}(1_{\tau})\Bigr).
\end{align*}
Also, define
\begin{align*}
\Phi_{\hat{G}}: C_{c}(\hat{G}) \to C_{b}(\hat{G}, M(\mathrm{fcph}(\mathcal{A}, G, \alpha)))
\end{align*}
by
\begin{align*}
\varphi \mapsto \Bigl(\tau \mapsto \varphi(\tau)\widetilde{\eta}_{\mathcal{L}}(1_{\tau})\Bigr).
\end{align*}
Finally, define
\begin{align*}
\Phi_{G}: C_{c}(G) \to C_{b}(\hat{G}, M(\mathrm{fcph}(\mathcal{A}, G, \alpha)))
\end{align*}
by
\begin{align*}
f \mapsto \Bigl(\tau \mapsto \widetilde{\eta}_{\mathcal{L}}(f_{\tau})\Bigr).
\end{align*}

\begin{figure}[h]
\begin{tikzcd}
	& {M(\mathcal{A}\rtimes\hat{G})} & {M(\mathcal{D})=M(\mathrm{fcph}(\mathcal{A}\rtimes\hat{G},G,\beta))} & G \\
	{\mathcal{A}} & {\hat{G}} \\
	&& {M(\mathcal{C})=M(\mathrm{fcph}(\mathcal{A},G,\alpha)\rtimes\hat{G})}
	\arrow[dashed, from=1-2, to=1-3]
	\arrow["{\Phi_{\mathcal{A}\rtimes \hat{G}}}", from=1-2, to=3-3]
	\arrow["{\Phi_{\mathcal{D}}}", dashed, from=1-3, to=3-3]
	\arrow[dashed, from=1-4, to=1-3]
	\arrow["{\Phi_G}", from=1-4, to=3-3]
	\arrow[dashed, from=2-1, to=1-2]
	\arrow["{\Phi_\mathcal{A}}"', from=2-1, to=3-3]
	\arrow[dashed, from=2-2, to=1-2]
	\arrow["{\Phi_{\hat{G}}}", from=2-2, to=3-3]
\end{tikzcd}
\caption{commutative diagram of construction $\Phi_{\mathcal{D}}$}
\end{figure}

Furthermore, by using the universality of the crossed product algebra \(\mathcal{A} \rtimes \hat{G}\) as indicated in the diagram, we construct the mapping \(\Phi_{\mathcal{A} \rtimes \hat{G}}\). Finally, by once again applying the universality of the full crossed product host \(\mathrm{fcph}(\mathcal{A}\rtimes\hat{G},G,\beta)\), Proposition \ref{prop: universal property of fcph}, we construct \(\Phi_{\mathcal{D}}\). The map \(\Phi_{\mathcal{D}}\) thus obtained satisfies
\[
\begin{aligned}
\mathcal{D} \ni \eta_{\mathcal{B}}(A\otimes \varphi)\eta_{\mathcal{L}'}(f)
&\mapsto \left(\tau \mapsto \varphi(\tau)\,\eta_{\mathcal{A}}(A)\,\eta_{\mathcal{L}}(f_{\tau})\right)\\[1mm]
&\in C_{c}(\hat{G}, \mathrm{fcph}(\mathcal{A},G,\alpha))
\subset \mathcal{C}.
\end{aligned}
\]

\medskip

\textbf{(The map $\Phi_{\mathcal{D}}$ is a homomorphism)}

It is obvious from the construction.

\medskip

\textbf{(The map $\Phi_{\mathcal{D}}$ is bijective)}

We construct the following non-degenerate $*$-homomorphism in a similar manner:
\begin{align*}
    \Psi_{\mathcal{C}}: \mathcal{C} \to M(\mathcal{D}).
\end{align*}
Since both $\Phi_{\mathcal{D}}$ and $\Psi_{\mathcal{C}}$ are non-degenerate, they can be uniquely extended to homomorphisms on the multiplier algebras:
\begin{align*}
    \widetilde{\Phi_{\mathcal{D}}}: M(\mathcal{D}) \to M(\mathcal{C}),\\
    \widetilde{\Psi_{\mathcal{C}}}: M(\mathcal{C}) \to M(\mathcal{D}).
\end{align*}
Since these maps are mutually inverses, they establish an isomorphism $M(\mathcal{C}) \simeq M(\mathcal{D})$.

The injectivity of $\Phi_{\mathcal{D}}$ is immediate from the injectivity of $\widetilde{\Phi_{\mathcal{D}}}$.

Now we show that $\Phi_{\mathcal{D}}$ is surjective. For any $c \in \mathcal{C}$, there exists some $m \in M(\mathcal{D})$ such that
\begin{align*}
    \widetilde{\Phi_{\mathcal{D}}}(m) = c.
\end{align*}
Let $(e_{\lambda})_{\lambda}$ be an approximate unit in $\mathcal{D}$. Then, we have $me_{\lambda} \to m$ strictly in $\mathcal{D}$, which implies
\begin{align*}
    \Phi_{\mathcal{D}}(me_{\lambda}) \to \widetilde{\Phi_{\mathcal{D}}}(m) = c, \quad \text{strictly}.
\end{align*}
Since the strict topology coincides with the norm topology on $\mathcal{C}$, we conclude that $c \in \mathrm{Im} \Phi_{\mathcal{D}}$, proving the surjectivity of $\Phi_{\mathcal{D}}$.

\end{proof}

\subsection{Takai duality in the case of a crossed product host}

As shown below, an isomorphism of the form $\mathcal{A} \otimes \mathcal{K}(L^2(G))$, similar to Takai duality, does not generally hold when the action is not continuous. However, a similar isomorphism can still be obtained.

We also assume that $G$ is a locally compact abelian group in this subsection.

Firstly, we define a $G$-product and Landstad's conditions.

\begin{dfn}[A $G$-product]
  Let $G$ be a locally compact abelian group. We say a C*-algebra $B$ is a \textbf{$G$-product} if the following two conditions are satisfied.
  
  (i) There is a homomorphism $\lambda: G\to UM(B)$ such that each function $g \mapsto \lambda_{g}y$, $y\in B$ is continuous from $G$ to $B$. 

  (ii) There is a homomorphism $\hat{\alpha}: \hat{G} \to \mathrm{Aut}(B)$ such that $(B,\hat{G}, \hat{\alpha})$ is a C*-dynamical system and 
  \begin{align*}
    \hat{\alpha}_{\tau}(\lambda_{g}) = \langle g, \tau\rangle \lambda_g, \quad \text{for all } g \in G, \tau \in \hat{G}.
  \end{align*}

  When given a $G$-product $\mathcal{B}$, we say that an element $x$ in $M(\mathcal{B})$ satisfies \textbf{Landstad's condtions} if:

  (L1) For all $\tau\in \hat{G}$, we have $\hat{\alpha}_\tau(x)=x$.

  (L2) For every $f$ in $L^1(G)$, we have $x\lambda_f \in \mathcal{B}$ and $\lambda_f x\in \mathcal{B}$.

  (L3) The map $g \mapsto \lambda_g x\lambda_{-g}$ is continuous on $G$.
\end{dfn}

\begin{rem}
  Let $(\mathcal{A}, G, \alpha)$ be a C*-action. Then, all crossed product hosts of $(\mathcal{A}, G, \alpha)$ are $G$-products.
\end{rem}

As it is known, a $G$-product can be expressed as a crossed product.

\begin{thm}[{\cite[Theorem~7.8.8]{MR3839621}}]
  A C*-algebra $\mathcal{B}$ is a $G$-product for a given abelian group $G$ if and only if there is a C*-dynamical system $(\mathcal{A}, G, \alpha)$ such that $\mathcal{B}=\mathcal{A} \rtimes_{\alpha} G$. This system is unique (up to isomorphism), and $\mathcal{A}$ consists of the elements in $M(\mathcal{B})$ that satisfy Landstad's conditions, whereas $\alpha_{g}=\lambda_{g}\cdot \lambda_{-g}, g\in G$. 
\end{thm}

For any crossed product host $\mathcal{C}$ of $(\mathcal{A}, G, \alpha)$, it can be easily checked that it is a $G$-product. By applying the above theorem, $C$ can be expressed as a form of crossed product,
\begin{align}
  \mathcal{C} \simeq \mathcal{A}^\# \rtimes_{\alpha^\#} G.
\end{align}

We call this C*-subalgebra $\mathcal{A}^\#$ of $M(\mathcal{C})$ \emph{Landstad algebra of $\mathcal{C}$}.

In \cite{MR3177335} a similar idea is used to define a subalgebra of $\mathcal{A}$, denoted by \(\mathcal{A}_0^{(\mathcal{L})}\), and the following relation holds:
\begin{align*}
  \mathcal{A}^\# \cap \eta_{\mathcal{A}}(\mathcal{A}) = \eta_{\mathcal{A}}(\mathcal{A}_0^{(\mathcal{L})}).
\end{align*}

\begin{thm}[Takai duality of C*-action]
  Let $(\mathcal{A},G,\alpha)$ be a C*-action, with $G$ abelian. Assume that there exists the full crossed product host $\mathrm{fcph}(\mathcal{A}, G, \alpha)$. Then, there is an isomorphism $\phi$ of $\mathrm{fcph}(\mathcal{A},G,\alpha)\rtimes_{\hat{\alpha}}\hat{G}$ onto $\mathcal{A}^\# \otimes K(L^2(G))$ such that the second dual action $\hat{\hat{\alpha}}$ of $G=\hat{\hat{G}}$ is carried into $\alpha^\#\otimes \mathrm{Ad} \rho$, where $\rho$ is the right regular representation of $G$ on $L^2(G)$.
\end{thm}
\begin{proof}
  By using the Landstad algebra $\mathcal{A}^\#$,
  \begin{align}
    \mathrm{fcph}(\mathcal{A},G,\alpha)\rtimes_{\hat{\alpha}}\hat{G} \simeq \mathcal{A}^\# \rtimes_{\alpha^\#} G \rtimes_{\hat{\alpha}} \hat{G}.
  \end{align} 
  Since $\hat{\alpha}=\widehat{\alpha^\#}$ on $\mathrm{fcph}(\mathcal{A},G,\alpha)=\mathcal{A}^\# \rtimes_{\alpha^\#} G$, we apply traditional Takai duality and get
  \begin{align}
    \mathcal{A}^\# \rtimes_{\alpha^\#} G \rtimes_{\widehat{\alpha^\#}} \hat{G} \simeq \mathcal{A}^\# \otimes K(L^2(G)).
  \end{align}
  These imply the desired isomorphism.
\end{proof}

\begin{rem}
  If $\alpha$ is not continuous, then the C*-action
\begin{align*}
  (\mathcal{A} \otimes K(L^2(G)),\, G,\, \alpha \otimes \operatorname{Ad}\rho)
\end{align*}
fails to be continuous. Together with the preceding results, this implies that, in general,
\begin{align*}
  \mathrm{fcph}(\mathcal{A}, G, \alpha) \rtimes_{\hat{\alpha}} \hat{G} \not\simeq \mathcal{A} \otimes K(L^2(G)).
\end{align*}
\end{rem}

In this way, it is quite interesting that taking the iterated crossed product host yields an induced continuous action of $G$.

In analogy with the case of crossed product algebras, we now consider when a given C*-algebra can be represented as a full crossed product host. To this end, we first introduce the following \emph{weak Landstad condition}.

\begin{dfn}[A weak Landstad Algebra]
Let $\mathcal{B}$ be a $G$-product. For $x\in M(\mathcal{B})$, consider the following weak form of the Landstad conditions:

  (L1) For all $\tau\in \hat{G}$, we have $\hat{\alpha}_\tau(x)=x$.

  (L2) For every $f$ in $L^1(G)$, we have $x\lambda_f \in \mathcal{B}$ and $\lambda_f x\in \mathcal{B}$.

  (L3') The map $g\mapsto \lambda_{g}x\lambda_{-g}$ is strictly continuous on $G$.

The set of all elements in $M(\mathcal{B})$ satisfying these conditions is called a \emph{weak Landstad algebra}. In particular, when $\mathcal{B}=\mathrm{fcph}(\mathcal{A},G,\alpha)$, we write $\mathcal{A}^b$ for the corresponding weak Landstad algebra.
\end{dfn}

Next, we describe some properties of the C*-algebra $\mathcal{A}^b$.

\begin{prop}
Let $\mathcal{B}=\mathrm{fcph}(\mathcal{A},G,\alpha)$. Then,
\begin{align*}
  \mathcal{A}^\# &\subset \mathcal{A}^b,\\
  \eta_{\mathcal{A}}(\mathcal{A}) &\subset \mathcal{A}^b.
\end{align*}
\end{prop}

\begin{proof}
Since condition (L3') is weaker than (L3), it immediately follows that
\begin{align*}
  \mathcal{A}^\# \subset \mathcal{A}^b.
\end{align*}
Next, by \cite[Corollary~8.4]{MR3177335}, the map
\begin{align*}
  G \to M(\mathcal{B}),\quad g \mapsto \eta_{\mathcal{A}}\bigl(\alpha_{g}(A)\bigr)
\end{align*}
is strictly continuous. Therefore, we obtain
\begin{align*}
  \eta_{\mathcal{A}}(\mathcal{A}) \subset \mathcal{A}^b.
\end{align*}
\end{proof}

It is thus expected that $\mathcal{A}^b$ is the largest C*-algebra which realizes $\mathrm{fcph}(\mathcal{A},G,\alpha)$ as a full crossed product host.

In the case of $G$ is non-abelian, we can have parallel considerations using the coaction.

By using \cite[Theorem~3.2]{MR2288078}, one obtains that for a C*-action $(\mathcal{A}, G, \alpha)$ with a full crossed product host $\mathrm{fcph}(\mathcal{A}, G, \alpha)$ with the maximal coaction, there exists a C*-algebra $\mathcal{A}^\#$ such that the following crossed product host representation holds:
\begin{align*}
  \mathrm{fcph}(\mathcal{A}, G, \alpha) \simeq \mathcal{A}^\# \rtimes_{\alpha^\#} G.
\end{align*}
In the commutative case, this $\mathcal{A}^\#$ coincides with the Landstad algebra, and hence the same notation is used.

Furthermore, we have the following:

\begin{thm}[Imai-Takai Duality for a full crossed product host]
Assume that a C*-action $(\mathcal{A}, G, \alpha)$ has the full crossed product host, where $G$ is not necessarily abelian. Then, the following isomorphisms hold:
\begin{align*}
  \mathrm{fcph}(\mathcal{A}, G, \alpha)\rtimes_{\hat{\alpha}} G &\simeq \mathcal{A}^\# \rtimes_{\alpha^\#} G \rtimes_{\hat{\alpha}} G\\
  &\simeq \mathcal{A}^\# \otimes K(L^2(G)).
\end{align*}
\end{thm}

\subsection{Example: The case $G=\mathbb{R}$}

Furthermore, we investigate the properties of $\mathcal{A}^\#$ in a concrete example.

\begin{example}
We now consider the example that was also used in \cite[Example~5.11]{MR3177335}, as follows. Let $\mathcal{H}$ be an infinite-dimensional Hilbert space, and set $\mathcal{A} := B(\mathcal{H}), G := \mathbb{R}$, and $\mathcal{L} := C^*(\mathbb{R})$. Let $H$ be a (possibly unbounded) self-adjoint operator, and define the action by
\begin{align*}
  \alpha_{t}(A) := U_{t} A U_{-t}, \quad t\in \mathbb{R},
\end{align*}
where \(U_{t} := e^{itH}\). A necessary and sufficient condition for the representation \((\pi = \operatorname{id}, U)\) of the C*-action \((\mathcal{A}, G, \alpha)\) to be cross is that \((i\mathbf{1} - H)^{-1} \in K(\mathcal{H})\). In this case, we have \(\eta_{\mathcal{L}}(\mathcal{L}) \subset K(\mathcal{H})\) and
\begin{align*}
  \mathcal{C}:=\mathrm{cph}(\mathcal{A}, G, \alpha) = K(\mathcal{H}).
\end{align*}
\end{example}

By definition,  
\begin{align*}
  \mathcal{A}^\# = \{X \in M(\mathcal{C}) = B(\mathcal{H}) \mid \text{\(X\) satisfies (L1), (L2), (L3)}\}.
\end{align*}

The right hand side is a *-subalgebra of $\mathcal{A}^\#$.

For example, it holds that
\begin{align*}
  \mathcal{A}^\# \supset \{X \in B(\mathcal{H}) \mid [X, H] \in B(\mathcal{H})\}.
\end{align*}
First, we show that any $X\in B(\mathcal{H})$ such that $[X,H]\in B(\mathcal{H)}$ belongs to \(\mathcal{A}^\#\), i.e., satisfies conditions (L1)-(L3).

(L1) Since the dual action \(\hat{\alpha}\) acts trivially on \(B(\mathcal{H})\), the condition (L1) is clearly satisfied.

(L2) Because both \(X\,\eta_{\mathcal{L}}(L)\) and \(\eta_{\mathcal{L}}(L)\,X\) lie in \(K(\mathcal{H}) = \mathcal{C}\), the condition (L2) is fulfilled.

(L3) If \([X, H] \in B(\mathcal{H})\), then the function
\begin{align*}
  X(t) = e^{itH}\,X\,e^{-itH}
\end{align*}
is differentiable with respect to \(t\) and hence continuous.

Therefore, any such \(X\) belongs to \(\mathcal{A}^\#\).

\subsection{Existence of ground states in C*-actions}

In this subsection, we discuss the existence of ground states in the setting of C*-actions.

The existence of ground states has been extensively studied. For instance, in the case of C*-dynamical systems on UHF algebras, it is well-known that the approximate innerness of the action provides a sufficient condition for the existence of a ground state (see, e.g., \cite{MR359623}). Moreover, analogous results have been obtained for more general C*-dynamical systems.

However, in certain settings such as the algebraic quantum field theory, one often encounters C*-actions that are not continuous, for example in the case of the Weyl C*-algebra. In these situations, it is natural to require analogous results within the framework of C*-actions, despite the lack of continuity. Although, in practice, one typically constructs representations by first obtaining a state satisfying the spectrum condition and then employing it to induce a representation, further investigations along this line may prove useful in broader contexts.

\medskip

In what follows, we aim to explore conditions under which ground states exist for such C*-actions, and to clarify how these conditions relate to the underlying algebraic structure, independently of the original C*-algebra.

In this work, we generalize the results of \cite{MR712636} to characterize the existence conditions for ground states in terms of ideals of the crossed product host.

\begin{thm}[Ground state conditions]\label{thm: gs}
  Let $(\mathcal{A}, \mathbb{R}, \alpha)$ be a C*-action, where $\mathcal{A}$ is a unital C*-algebra. And assume that there exists the full crossed product host $\mathrm{fcph}(\mathcal{A}, \mathbb{R}, \alpha)$. Then the following conditions are equivalent:

  (i) There is an $\alpha$-invariant ground state of $\mathcal{A}$.

  (ii) There is a proper, closed, two-sided ideal $I$ of $\mathrm{fcph}(\mathcal{A}, G, \alpha)$ such that $I \subset \hat{\alpha}_{\lambda}(I)$ for any $\lambda\geq 0$ and the union of $\hat{\alpha}_\lambda(I)$ with all $\lambda\in \mathbb{R}$ is dense in $\mathrm{fcph}(\mathcal{A}, G, \alpha)$.
\end{thm}
\begin{proof}
  In \cite[Theorem~3.5]{MR712636}, the crossed product $\mathcal{A} \rtimes \mathbb{R}$ is used; however, here we consider $\mathrm{fcph}(\mathcal{A}, \mathbb{R}, \alpha)$ as the crossed product host and carry out an analogous construction. Finally, we obtain the desired state by constructing a state on $B(\mathcal{H})$ as a limit and then pulling it back to a state on $\mathcal{A}$.
  
  The spectral evaluation is carried out in a similar manner.
\end{proof}

\subsection{A reduced crossed product host}

In this subsection, we consider a reduced crossed product host, which serves as the counterpart of the reduced crossed product in the context of crossed product hosts.

Let $(\mathcal{A}, G, \alpha)$ be a C*-action, and let $(\pi, \mathcal{H})$ be a non-degenerate representation of $\mathcal{A}$. Then, the covariant representation $(\widetilde{\pi}, \lambda)$ of the C*-action on $L^2(G,\mathcal{H})$ is given by
\begin{align*}
    (\widetilde{\pi}(A)\xi)(g) &:= \pi(\alpha_{g^{-1}}(A))\xi(g),\\
    (\lambda_{h}\xi)(g) &:= \xi(h^{-1}g),
\end{align*}
where $\xi\in L^2(G,\mathcal{H})$.

\begin{dfn}
Let $(\pi_u, \mathcal{H}_{u})$ be the universal representation of $\mathcal{A}$. When the covariant representation $(\widetilde{\pi}_{u}, \lambda)$ of the C*-action $(\mathcal{A}, G, \alpha)$ is cross, we call the crossed product host constructed from $(\widetilde{\pi}_{u}, \lambda)$ the \textbf{reduced crossed product host}, denoted by $\mathrm{rcph}(\mathcal{A},G,\alpha)$. That is,
\begin{align*}
    \mathrm{rcph}(\mathcal{A},G,\alpha):=C^*(\widetilde{\pi}_{u}(\mathcal{A})\lambda_{\mathcal{L}}(\mathcal{L})).
\end{align*}
\end{dfn}

\begin{lemma}\label{lem: amenable}
Suppose that the C*-action $(\mathcal{A},G,\alpha)$ has a full crossed product host. If the group $G$ is amenable, then for any covariant representation $(\pi, U)$, we have
\begin{align*}
    \|\pi(A)U_{\mathcal{L}}(L)\| \leq  \|\widetilde{\pi}_{u}(A)\lambda_{\mathcal{L}}(L)\| \quad \text{for all } A\in \mathcal{A}, L\in \mathcal{L}.
\end{align*}
\end{lemma}

\begin{proof}
  Let $z \in B(L^2(G,\mathcal{H}))$ be a unitary operator satisfying $(z \xi)(g)=U(g)^{-1}(\xi(g))$. Then, we have
  \begin{align*}
    z(\lambda_{g} \otimes U_{g})z^{-1} &= \lambda_{g} \otimes 1,\\
    z(1 \otimes \pi(A)) z^{-1} &= \pi(A),
  \end{align*}
  which leads to
  \begin{align*}
    \|\widetilde{\pi}(A)\lambda(f)\| = \|(1 \otimes \pi)(A)(\lambda \otimes U)(f) \|.
  \end{align*}
  
  Let $\varepsilon>0$ be arbitrary. We aim to show that
  \begin{align*}
    \|(1 \otimes \pi)(A)(\lambda \otimes U)(f) \| \geq \|\pi(A)U_{\mathcal{L}}(f)\| - \varepsilon  \quad \text{for all } A\in \mathcal{A}, f\in C_{c}(G).
  \end{align*}
  
  First, choose $\xi_{0} \in \mathcal{H}$ such that
  \begin{align*}
    \|\xi_{0}\|=1, \quad \|\pi(A)U_{\mathcal{L}}(f)\xi_{0}\| > \|\pi(A)U_{\mathcal{L}}(f)\| - \frac{\varepsilon}{2}.
  \end{align*} 
  Define
  \begin{align*}
    \delta = \left(\frac{\|\pi(A)U_{\mathcal{L}}(f)\| - \frac{\varepsilon}{2}}{\|\pi(A)U_{\mathcal{L}}(f)\| - \varepsilon}\right)^2 -1.
  \end{align*}
  Clearly, $\delta>0$. Define the subset $S \subset G$ as $S:=\mathrm{supp}(f) \cup \{1\}$, which is compact. Since $G$ is amenable, we can use the F\o lner condition. That is, there exists a compact subset $K \subset G$ such that
  \begin{align*}
    0 < \mu(K) < \infty, \quad \frac{\mu(S^{-1}K \triangle K)}{\mu(K)} < \delta.
  \end{align*}
  Since $1 \in S^{-1}$, we obtain $\mu(S^{-1}K \setminus K) < \delta \mu(K)$, which implies $\mu(S^{-1}K) < (1+\delta) \mu(K)$.
  
  Next, define $\xi\in L^2(G,\mathcal{H})$ as
  \begin{align*}
    \xi(g) = \begin{cases}
    \xi_{0}, & g\in S^{-1}K,\\
    0, & g\notin S^{-1}K.
    \end{cases}
  \end{align*}
  Then,
  \begin{align}\label{eq: xi norm}
    \|\xi\| = \mu(S^{-1}K)^{1/2}\|\xi_{0}\| < (1+\delta)^{1/2}\mu(K)^{1/2}.
  \end{align}
  Using this, we estimate $\|(1 \otimes \pi)(A)(\lambda \otimes U)(f) \xi\|$.
  
  For $g\in K$, we have
  \begin{align*}
    ((1 \otimes \pi)(A)(\lambda \otimes U)(f) \xi)(g) &= \int_{G}((1 \otimes \pi)(Af(h)) (\lambda_{h} \otimes U_{h})\xi)(g) d h\\
    &= \int_{G}\pi(Af(h))U_{h}\xi(h^{-1}g)d h\\
    &= \int_{G}\pi(A)f(h)U_{h}\xi_{0}d h\\
    &= \pi(A) U(f)\xi_{0}.
  \end{align*}
  Here, in the third equality, we used the fact that if $h \in \mathrm{supp}(f) = S$, then $\xi(h^{-1}g) = \xi_0$.
  
  Thus,
  \begin{align*}
    \|(1 \otimes \pi)(A)(\lambda \otimes U)(f) \xi\| &\geq \mu(K)^{1/2}\|\pi(A) U(f)\xi_{0}\| \\
    &>  \mu(K)^{1/2}\left(\|\pi(A)U_{\mathcal{L}}(f)\| - \frac{\varepsilon}{2}\right).
  \end{align*}
  From equation \eqref{eq: xi norm}, we obtain
  \begin{align*}
    \|(1 \otimes \pi)(A)(\lambda \otimes U)(f) \xi\| &> \frac{\mu(K)^{1/2}(\|\pi(A)U_{\mathcal{L}}(f)\| - \frac{\varepsilon}{2})}{(1+\delta)^{1/2} \mu(K)^{1/2}}\\
    &= (1+ \delta)^{-1/2}\left(\|\pi(A)U_{\mathcal{L}}(f)\| - \frac{\varepsilon}{2}\right)\\
    &= \|\pi(A)U_{\mathcal{L}}(f)\| - \varepsilon.
  \end{align*}
  This completes the proof.
\end{proof}

\begin{thm}
Suppose that the C*-action $(\mathcal{A}, G, \alpha)$ has a full crossed product host. If the group $G$ is amenable, then the full crossed product host and the reduced crossed product host coincide. That is,
\begin{align*}
    \mathrm{fcph}(\mathcal{A},G,\alpha) \simeq \mathrm{rcph}(\mathcal{A},G,\alpha).
\end{align*}
\end{thm}

\begin{proof}
By \cite[Theorem~5.8]{MR3177335}, there exists a unique $*$-homomorphism
\begin{align*}
    \Phi: \mathrm{fcph}(\mathcal{A},G,\alpha) \to \mathrm{rcph}(\mathcal{A},G,\alpha)
\end{align*}
such that
\begin{align*}
    (\mathrm{rcph}(\mathcal{A},G,\alpha), \eta_{\mathcal{A},r}, \eta_{\mathcal{L}, r}) = (\Phi(\mathrm{fcph}(\mathcal{A},G,\alpha)), \widetilde{\Phi}\circ\eta_{\mathcal{A}},\widetilde{\Phi}\circ\eta_{\mathcal{L}}).
\end{align*}
This implies that $\mathrm{rcph}(\mathcal{A},G,\alpha)$ is a quotient C*-algebra of $\mathrm{fcph}(\mathcal{A},G,\alpha)$.

Since $\Phi$ is surjective, it suffices to show that it is isometric.
By the Lemma \ref{lem: amenable}, for any covariant representation $(\pi,U)$, we have
\begin{align*}
    \|\pi(A)U_{\mathcal{L}}(L)\| \leq  \|\widetilde{\pi}_{u}(A)\lambda_{\mathcal{L}}(L)\|.
\end{align*}
Since $\mathrm{span}\{\pi(\mathcal{A})U_{\mathcal{L}}(\mathcal{L})\}$ is dense in $\mathrm{rcph}(\mathcal{A},G,\alpha)$, and $\Phi$ is continuous, it follows that for any $C\in \mathrm{rcph}(\mathcal{A},G,\alpha)$,
\begin{align*}
    \|\Phi(C)\| \leq \|C\|.
\end{align*}
Since $\Phi$ is a $*$-homomorphism, it is norm-decreasing, hence
\begin{align*}
    \|C\| \leq \|\Phi(C)\|.
\end{align*}
Thus, $\|\Phi(C)\| = \|C\|$, proving that $\Phi$ is an isometry and hence an isomorphism.
\end{proof}
\section{Conclusion}

As we have seen throughout this paper, the crossed product host $\mathcal{C}$ possesses the structure of a \( G \)-product and can thus be expressed as a crossed product algebra $\mathcal{A}^\# \rtimes_{\alpha^\#} G$ via the Landstad algebra $\mathcal{A}^\#$. Utilizing this, as demonstrated in Theorem \ref{thm: gs}, it is expected that proofs independent of the C*-algebra $\mathcal{A}$ will also hold for the crossed product host.

In this study, we have focused primarily on the cases where the host algebra is given by $\mathcal{L} = C^*(G)$ or where a full crossed product host exists. However, natural extensions to more general settings can be considered.

Furthermore, one of our initial objectives was to establish a suitable condition ensuring the existence of a ground state for a C*-action. Although we have provided a necessary and sufficient condition, it is not particularly convenient for practical applications. Therefore, further research is needed to formulate a more usable criterion, especially for C*-actions such as those on the Weyl algebra.

Additionally, the representation of a given C*-algebra as the full crossed product host of a C*-action is not necessarily unique. The extent to which such C*-actions exist, as well as the relationships between different C*-actions, remain largely unexplored. Investigating these aspects will be an important step toward a deeper understanding of the structure and classification of crossed product hosts.
\section*{Acknowledgment}

I would like to sincerely thank my supervisor, Professor Yasuyuki Kawahigashi, for his support, encouragement, and invaluable guidance throughout this research. It was supported by the WINGS-FMSP program at The University of Tokyo.
\bibliographystyle{alpha}
\bibliography{reference}

\begin{thebibliography}{EKQR06}

\bibitem[EKQR06]{MR2203930}
Siegfried Echterhoff, S.~Kaliszewski, John Quigg, and Iain Raeburn.
\newblock A categorical approach to imprimitivity theorems for {$C^*$}-dynamical systems.
\newblock {\em Mem. Amer. Math. Soc.}, 180(850):viii+169, 2006.

\bibitem[GN14]{MR3177335}
Hendrik Grundling and Karl-Hermann Neeb.
\newblock Crossed products of {$C^*$}-algebras for singular actions.
\newblock {\em J. Funct. Anal.}, 266(8):5199--5269, 2014.

\bibitem[GN20]{MR4157367}
Hendrik Grundling and Karl-Hermann Neeb.
\newblock Crossed products of {$C^*$}-algebras for singular actions with spectrum conditions.
\newblock {\em J. Operator Theory}, 84(2):369--451, 2020.

\bibitem[IT78]{MR500719}
Sh\=o{} Imai and Hiroshi Takai.
\newblock On a duality for {$C\sp{\ast} $}-crossed products by a locally compact group.
\newblock {\em J. Math. Soc. Japan}, 30(3):495--504, 1978.

\bibitem[KQ07]{MR2288078}
S.~Kaliszewski and John Quigg.
\newblock Landstad's characterization for full crossed products.
\newblock {\em New York J. Math.}, 13:1--10, 2007.

\bibitem[Kus83]{MR712636}
Masaharu Kusuda.
\newblock Crossed products of {$C\sp{\ast} $}-dynamical systems with ground states.
\newblock {\em Proc. Amer. Math. Soc.}, 89(2):273--278, 1983.

\bibitem[Ped18]{MR3839621}
Gert~K. Pedersen.
\newblock {\em {$C^*$}-algebras and their automorphism groups}.
\newblock Pure and Applied Mathematics (Amsterdam). Academic Press, London, second edition, 2018.
\newblock Edited and with a preface by S\o ren Eilers and Dorte Olesen.

\bibitem[PS75]{MR359623}
Robert~T. Powers and Sh\^oichir\^o Sakai.
\newblock Existence of ground states and {KMS} states for approximately inner dynamics.
\newblock {\em Comm. Math. Phys.}, 39:273--288, 1974/75.

\bibitem[Rae88]{MR989764}
Iain Raeburn.
\newblock On crossed products and {T}akai duality.
\newblock {\em Proc. Edinburgh Math. Soc. (2)}, 31(2):321--330, 1988.

\bibitem[Tak75]{MR365160}
Hiroshi Takai.
\newblock On a duality for crossed products of {$C\sp{\ast} $}-algebras.
\newblock {\em J. Functional Analysis}, 19:25--39, 1975.

\bibitem[Wil07]{MR2288954}
Dana~P. Williams.
\newblock {\em Crossed products of {$C{^\ast}$}-algebras}, volume 134 of {\em Mathematical Surveys and Monographs}.
\newblock American Mathematical Society, Providence, RI, 2007.

\end{thebibliography}

\end{document}